\documentclass[]{interact}

\usepackage{listings}
\usepackage{xcolor}
\usepackage{inconsolata} 
\definecolor{darkpink}{RGB}{255, 153, 204}
\definecolor{lightcyan}{RGB}{0, 180, 180}
\definecolor{darkred}{RGB}{0.8, 0, 0}
\definecolor{pykeyword}{RGB}{0,85,160}
\definecolor{pycomment}{RGB}{106,106,106}
\definecolor{pystring}{RGB}{170,0,85}
\definecolor{pynumber}{RGB}{128,0,128}
\definecolor{pylight}{RGB}{248,248,255}
\definecolor{smoothgreen}{RGB}{0,153,0}
\lstdefinestyle{PythonCustom}{
language=Python,
backgroundcolor=\color{pylight},
basicstyle=\ttfamily\footnotesize,  
keywordstyle=\color{pykeyword}\bfseries,
commentstyle=\color{pycomment}\itshape,
stringstyle=\color{pystring},
numberstyle=\tiny\color{pycomment},
numbers=left,
numbersep=6pt,                
frame=none,
framesep=4pt,                  
rulecolor=\color{pycomment},
breaklines=true,
breakatwhitespace=false,
showstringspaces=false,
tabsize=4,
captionpos=none,
xleftmargin=0pt,
xrightmargin=0pt,
linewidth=\linewidth,
aboveskip=6pt,
belowskip=6pt,
escapeinside={(*@}{@*)},
morekeywords={self,True,False,None,async,await,with,as,from,import,def,class,return,if,elif,else,for,while,try,except,finally,raise,lambda,pass,break,continue,yield,assert,del,in,is,and,or,not,P,I,RI,R},
emph=[2]{0,1,2,3,4,5,6,7,8,9},
emphstyle=[2]\color{pynumber},
morestring=[b]',
morestring=[b]",
}
\usepackage{epstopdf}
\usepackage[caption=false]{subfig}

\usepackage[natbibapa,nodoi]{apacite}

\theoremstyle{plain}

\theoremstyle{definition}

\theoremstyle{remark}

\begin{document}


\title{Proliferating series by Jean Barraqué: a study and classification in mathematical terms}

\author{
	\name{Isabel Tardón\textsuperscript{a,b} and Pablo Martín-Santamaría\textsuperscript{a,b}\thanks{CONTACT Isabel Tardón, Email: isabel.tardon.b@gmail.com}}
	\affil{\textsuperscript{a}Universidad de Málaga, Malaga, Spain}
	 \affil{\textsuperscript{b}Conservatorio Superior de Música de Málaga, Malaga, Spain }
}

\maketitle

\begin{abstract}
	\label{abstract}
	
Barraqué's proliferating series give an interesting turn on the concept of classic serialism by creating a new invariant when it comes to constructing the series: rather than the intervals between consecutive notes, what remains unaltered during the construction of the proliferations of the given base series is the permutation of the notes which happens between two consecutive series, that is to say, the transformation of the order of the notes in the series. This presents new possibilities for composers interested in the serial method, given the fact that the variety of intervals obtained by this method is far greater than that of classic serialism. 

In this manuscript, we will study some unexplored possibilities that the proliferating series offer from a mathematical point of view, which will allow composers to gain much more familiarity with them and potentially result in the creation of pieces that take serialism to the next level.

\end{abstract}


\section{Introduction}
\label{sec:Introduction}

Since 1923, Schönberg’s dodecaphonic method has allowed composers to write non-tonal music in a methodic and coherent way. The Second Viennese School developed the earliest pieces using the dodecaphonic method, leading to works of utmost importance, such as the Variations for Orchestra Op. 31 by Schönberg, Alban Berg's Violin Concerto or the Concerto for Nine Instruments by Anton Webern. After World War II, members of the Darmstadt School continued using this method and expanded the concept of a series, which originally referred to the height of the notes, to cover other aspects such as dynamics, articulations, note lengths, etc (\cite{bandur2001aesthetics}). 

During several summer courses organized by Le Guide du Concert during the summers of 1956 to 1959, a new approach to serialism was exposed to several students by Jean Barraqué, as \cite{riotte1987series} tells, which is what we will address in this article. He suggested a brand new way to generate material by repeatedly applying the transformation which relates two series in order to generate other series with a wide range of intervallic structures. The result of this idea was the concept of proliferating series. 
In his work “... Au Delà du Hasard”, based on the literary text “La Mort du Virgil” (\cite{ozzard1989barraque}), he uses proliferating series to create various pieces, some of which have two generator series that are related by the usual transformations of classical serialism.

The question then arises as to how these related series behave when we apply the proliferating process: how many series can we generate? What is their structure? How can we use the answers to our previous questions to our advantage in order to create music?
In this article, we answer these questions, using mathematical group theory as a basis with which to explain the behaviour of different transformations in order to give composers a thorough knowledge of the properties of these series. 
In order to do this, we will talk about the way in which the series are constructed, the possible numbers of proliferations and structures we can obtain using two related series, how to build a series with the desired properties so that when it is matched with its retrograde inversion it produces however much material is needed and we will conclude with the relevance of these findings for composers and theoretical musicians, ending with a conclusion to sum up the main results.

\section{Methodology}
\label{sec:methodology}

As said in the introduction, we are going to analyse the structural properties of a certain type of proliferating series using mathematical group theory. This means that most of this article will use rather technical arguments for those not acquainted with mathematics. Our goal, though, is to expose those ideas as clearly as possible for musicians to understand them, so we will carefully explain every concept that we use and highlight its implications with proliferating series.

In particular, we will study the proliferating series constructed from the transformations of traditional serialism, so we will independently treat each transformation and each of their transpositions, grouping those that share common properties.

A key tool for our investigation was the Python script presented in Appendix~\ref{sec: appendix}. Besides a huge improvement in series manipulation, we have obtained empirical data to motivate and support the mathematical laws behind proliferation. We leave the code here as a tool for composers or other musicians to use proliferating series even without much idea of their construction, and to obtain the data that verifies our study. 

Most of the examples used to illustrate our reasonings about proliferating series will use remarkable series used in traditional serialism by known composers, showing how this process would have worked for those and what new material it brings to us.

\section{Construction of the proliferating series}
\label{sec:construction}
	
The construction of Barraqué's proliferating series, as exposed in \cite{nicolas1987souci}, consists of choosing two conventional dodecaphonic series and then taking the transformation of the notes required to go from the first to the second of the chosen series for the construction of the subsequent series. 

\begin{figure}[htb] 
	\centerline{
		\resizebox*{15cm}{!}{\includegraphics{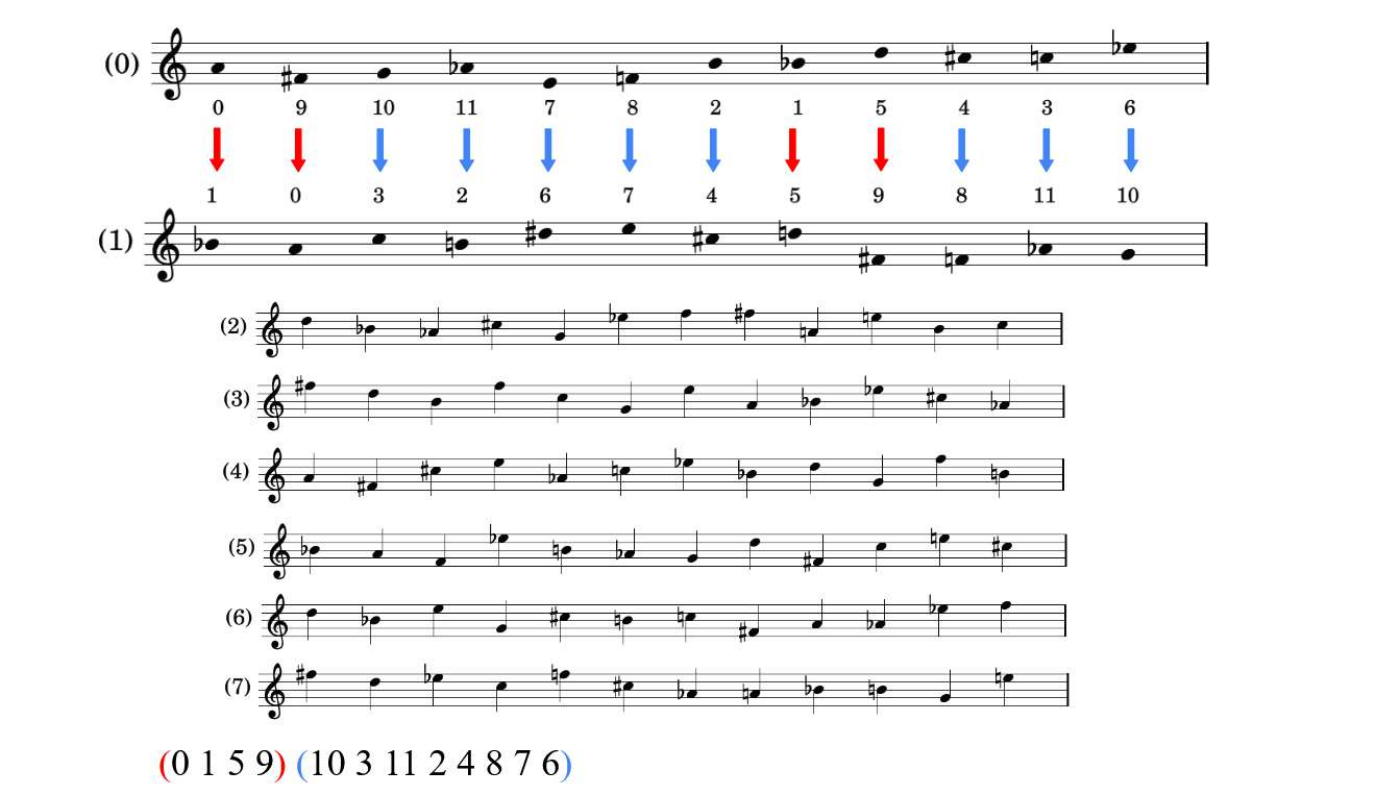}}\hspace{0pt}}
	\caption[Example of the proliferation process starting from two unrelated series and repeating until all the proliferations are obtained.]{Mapping the series from Webern Symphony Op. 21 and Webern String Quartet Op. 28 in order to obtain a permutation which generates proliferations}
	\label{fig: construction example}
	
\end{figure}

To illustrate this, let us take two arbitrary series, placing the second series below the first one, and consider the arrows which map one note of the first series to one on the second, as shown in Figure~\ref{fig: construction example}. 

The arrows in the figure define a permutation, that is to say, a function that changes the order of the twelve notes. In this example, both of the series used belong to Webern: $[A, F\#, G, Ab, E, F, B, Bb, D, C\#, C, Eb]$ is taken from his Symphony Op. 21, and $[Bb, A, C, B, D\#, E, C\#, D, F\#, F, Ab, G]$ is from his String Quartet Op. 28. The permutation which relates the two series transforms the note $A$ into $Bb$, $F\#$ into $A$, $G$ into $C$, etc.

If we apply the permutation to the first series, we obtain, by definition, the second series. But what happens when we apply this same permutation to the second series? We get a totally new third series, which has a $B$ where the second had a $A$, has an $A$ where the second had an $F\#$ and so on. We can apply our permutation again to the third series to obtain a fourth series, a fifth etc. We call those new series the proliferations of the initial pair. When repeating this process, we observe that eventually we get to the original series, as we can see in Figure~\ref{fig: construction example}; this will always happen due to elemental properties of group theory (\cite{dixon1996permutation}). The total number of series obtained with this method is what we call the order of the permutation. In the example presented in Figure~\ref{fig: construction example}, the order is 8.

For future reasoning, it is particularly useful to express the permutation we have obtained as a product of disjoint cycles. Each of the cycles is represented in Figure~\ref{fig: construction example} using arrows of a different colour. In a cycle, each note gets transformed into the one written immediately after it. For the last note of the cycle, the note that follows is the first of the same cycle. 
For example, let us look at the cycle $(0 \; 1\; 5\; 9)$ from Figure~\ref{fig: construction example}, which corresponds to the notes $(A\; Bb\; D\; F\#)$. It indicates the transformation that turns $A$ into $Bb$, $Bb$ into $D$, $D$ into $F\#$, and $F\#$ back into $A$. Expressing a permutation as a product of disjoint cycles consists of nothing more than writing each of these cycles one after the other, separated by parentheses, in no particular order. This can be done with any permutation and it is a clear depiction of which notes are changed in the permutation. 
Furthermore, the order of the permutation (which is the number of series that we obtain by proliferation) is the least common multiple of the lengths of the cycles involved in the permutation (\cite{dixon1996permutation}). 

We say that the structure of a permutation is the combination of cycles of a certain length that matches its representation in disjoint cycles, and we will write them as [\{$C_1$\}, \{$C_2$\}, …], where \{$C_i$\} represents an abstract cycle of length $C_i$. In Figure~\ref{fig: construction example}, we can see how each cycle is represented by arrows of a different colour. It has two cycles, of length 4 and  8, so we say that it has the structure [\{$4$\}, \{$8$\}].

\subsection{Series constructed through the transformations of traditional serialism}
\label{subsec: trad construction}

In traditional twelve-tone serialism, the series chosen as the main material for a composition can be transformed through any combination of transposing, inverting and retrograding said series. While in the construction method described above the two series used to generate all the others need not have any relation between them, it is of great mathematical as well as musical interest to cause a relation by using these transformations. 

The series that we obtain through these transformations are classified into four types, depending on the proceedings that need to be applied to the original series in order to obtain them. We list these types below along with the abbreviation that we will use in order to refer to them:

\begin{enumerate}
	\item[(P)] Primes: They are obtained by only applying a transposition $t$ to the original series.
	\item[(R)] Retrogrades: They are the result of reading the series in reverse order (back to front) and applying a transposition $t$.
	\item[(I)] Inversions: They are obtained by fixing the first note, then exchanging each interval in the original series for its complementary and finally applying a transposition $t$ to the resulting series.
	\item[(RI)] Retrograde Inversions: They are obtained by applying retrograding and inversion and then applying a transposition $t$ to the result.
	
\end{enumerate}

Note that $t=0$ is a transposition in all cases.

Therefore, rather than choosing two arbitrary unrelated series, we will pick one series and then choose the second from the list of primes (P), retrogrades (R), inversions (I) and retrograde inversions (RI) of the first. The construction which uses the retrograde inversion is mentioned in \cite{hopkins1978barraque} in relation to the piece "...Au Delà du Hasard". 

\begin{figure}[htb]  
	\centerline{
		\resizebox*{14cm}{!}{\includegraphics{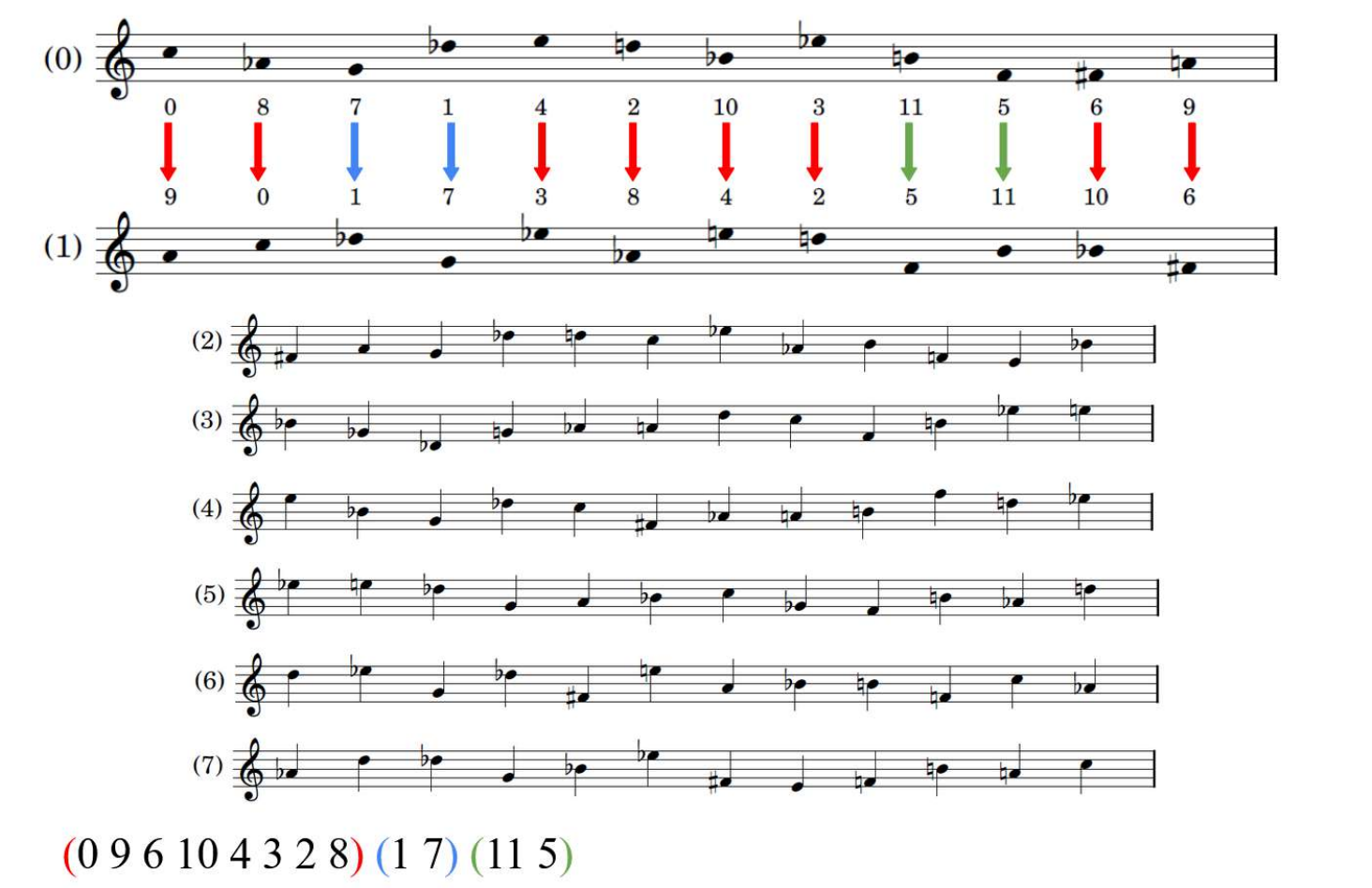}}\hspace{0pt}}
	\caption[Proliferation process using a series and its retrograde inversion to showcase the new material that results from combining traditional serialism with proliferating serialism.]{Proliferations of the series from "...Au Dela du Hasard", as used by Jean Barraqué}
	\label{fig: barraque construction}
	
\end{figure}

Aside from the relation between the series, the method of construction remains unaltered: once again, we take the two series which are as explained above and consider the transformation of the notes which occurs from the first to the second. As we can see in Figure~\ref{fig: barraque construction}, we get a number of series with a wide variety of intervallic structures, just as we did before in Figure~\ref{fig: construction example}.

This method of proliferation will be the main focus of this paper from this point forth, as we will study the number of proliferations that they can produce, the possible structures and the way in which a composer can build an initial series so that it fits their creative needs. To avoid confusion with other uses of permutations and for simplicity, we will call the permutations obtained from proliferation from a given series to one of its transformations through traditional serialism “proliferating permutations”(PPs).

When working with these series, we will label the notes in the usual way for dodecaphonic series: the first note of the first series shall be assigned the number zero and, for all other notes, their assigned number will be the distance to that first note in number of semitones (the series will be numbered modulo 12). 

Therefore, when we apply a transformation to a series, the result is as follows:

\begin{enumerate}
	\item[(i)] For retrograding, we simply rewrite the original string of numbers in reverse order.
	\item[(ii)] For transposition, we add to every number in the series the number of semitones which we wish to transpose and then write those numbers modulo 12, adding or subtracting 12 until the number is between 0 and 11.
	\item[(iii)] For inversion, rather than counting semitones in an ascending order, we will do so in descending order, so the set of pitches we obtain will have the first number be zero and for the rest of the numbers, we will write the same as in the original series, preceded by a negative sign. In order to go back to a numbering that is from 0 to 11, we simply add 12 to all negative numbers. Note that this works only because our series starts with 0. In any other case, the behaviour would be analogous, though the notation would become slightly more cumbersome.
	Another way to see this is that, for every number in the original series, the inverted series has the number such that, when it is added to the one in the original series, the result is 12 (in this context, 12 and 0 act as the same number).
\end{enumerate}

\subsection{Generalization of proliferating series: micro-tonal scales}
\label{subsec: microtonal}

By using modulo $n$ instead of modulo 12 and applying the same proceedings, we generalise the concepts explained above, extending the possibilities of the method to divisions of the octave that are different from the chromatic scale. These include microtonal divisions, such as quarter-tones or third-tones, modal scales such as the whole-tone scale, and overall any division imaginable. 
Now, instead of counting semitones, we will count the number of notes in our division between the first note of the series and the rest of the notes (the concept is the same as in the chromatic circle, except instead of having 12 notes, there can be any number of notes that we wish). For simplicity, we may refer to this microtonal unit of measurement as a tone-fraction. 

Essentially, what this means is that, given a set of n pitches, we can establish a “natural” order in the set whereby the distance between two given pitches is the number of pitches that we pass through in this natural order to go from the first to the second. Since this order is cyclical, when we write a series (a reordering of the pitches), we simply label the first pitch of the series as 0 and then place the remaining labels by counting the distance between pitches, always relating back to the natural order that we chose initially.

It is also worth mentioning that rhythms or any other musical parameter observed in integral serialism can receive the same treatment by numbering the musical elements in whichever way we see fit. The construction of proliferating series using these other divisions and parameters is analogous to the one mentioned before, in which the octave is divided into twelve semitones.

\section{Possible number of proliferations for a series of $n$ notes}
\label{sec: Prolif n notes}

When beginning a serial composition, it is crucial for the composer to know the material that they have at their disposition. In the case of traditional twelve-tone serialism, this is fairly straightforward: a 12 by 12 matrix is created where the composer can immediately see the series they have chosen as well as all possible transpositions, inversions and retrogrades. Different properties series used in this method have already been widely studied, first by the Second Viennese School and afterwards by the Darmstadt School and other musicians.

Going back to proliferating series, if one uses two arbitrary unrelated series with $n$ notes, the only limit to the number of proliferations that can be obtained is the orders that the elements of the symmetric group of size $n$ (the group of all the permutations of $n$ elements) have. However, the construction which uses a series and one of its transformations derived from classical serialism has further limitations regarding the number of proliferations due to the fact that the two generating series are not independent from one another. That is, only some of the permutations of the symmetric group can be a PP.

Thus, certain questions arise: 
Do the PPs obtainable by each transformation all follow a common structure?
More specifically, how many proliferations can be obtained using a series of $n$ notes and one of its transformations? 
We will answer these questions from a systematic point of view, using the tools described before. In each of the 4 subsections we will examine the transposite variants of each transformation, using t to express the number of transpositions and $n$ for the number of notes.

\subsection{Study of proliferations using inversion (I)}
\label{sec: Prolif I}

This first case is the simplest one. When we invert a series and then transpose it, we will get a PP of order (at most) two, no matter the number of transpositions, meaning that, in the proliferating process, we will obtain the original series again when we apply the PP to the inversion. 

To prove this statement, let $t$ be the transposition used and $x$ be the number corresponding to one arbitrary note. The PP will transform $x$ into $-x+t$, again, because of first applying an inversion to a series that starts with $0$ and then adding a transposition. Notice that, in turn, this note $-x+t$ will be transformed into $-(-x+t)+t=x$, and that means that $(x, -x+t)$ is a cycle (of length two) of the PP. Also notice that if we first transpose our series and then invert it, the resulting series is the same; we start from $x$, we transpose it to $x+t$, and to make the inversion, since our series now starts with $t$, we get $2t-(x+t)=-x+t$, the same as in the other case. In mathematical terms, we say that inversion and transposition commute. We picked an arbitrary number of transpositions and an arbitrary note, and we got a cycle of order two. This means that those will be the only possible cycles of a PP obtained like this, not including the degenerated case in which $x=-x+t$, which would make the cycle of length 1. We get that the PP will be the product of cycles of length 1 and 2, so its overall order will also be 1 or 2, being 1 just in the case of a series of one note or two notes with 0 transpositions.

As we can see, this case brings us no new material whatsoever.

\subsection{Study of proliferations using primes (P)}
\label{sec: Prolif P}

This case is also rather simple. Of course, using P0 (the unaltered series from which we derive the rest) will generate a PP that is the identity permutation, which leaves all the notes as they are and provides no new material. But in any other case, the PP transforms an arbitrary note $x$ into $x+t$. This will happen for all the notes, so if we apply the PP twice, we will get $(x+t)+t=x+2t$. In general, when we apply the PP an arbitrary number $k$ times, we get $x+kt$. And to find how many times we have to apply it to return to the note $x$ and finish the cycle, we just have to determine when $x=x+kt$ modulo $n$. 
This will happen after $k=n/gcd(n,t)$ iterations. That is to say, just as in the previous case (with the exception of the few cycles of length 1, which won't happen in this case), every cycle of our PP will have the same length, and the order of the permutation is precisely that length: $n/gcd(n,t)$. The proliferations will be just transpositions of the original, so, again, it does not bring new material.

\subsection{Study of proliferations using retrograde inversion (RI)}
\label{sec: Prolif RI}

For this section, we must first discuss in which order to apply retrogradation, inversion and transposition. We already proved that inversion and transposition commute, and it is quite evident to observe that retrogradation and transposition also commute. By contrast, inversion and retrogradation do not: note that if the process to be followed is first the inversion of the series and then its retrograding, when we start with an initial series whose first note is assigned zero, the inverted and then retrograded version will have the number zero in its last position. 
If, on the contrary, the retrograding of the series is done before the inversion, the result is a series that presents the same interval distribution as before, but with a different transposition.

As it turns out, the structural properties are exactly the same using one order or the other. Just for this reasoning, we will call RI the process of first inverting and then retrograding and IR the inverse order. Observe that applying RI and a transposition $t$ and then IR and a transposition $n-t$ to a series yields the same series. If we consider the series of $n=7$ notes $A=[0, 3, 4, 2, 1, 6, 5]$, inverting and then retrograding, with a transposition $t=2$, yields the series $B=[4, 3, 1, 0, 5, 6, 2]$. If we now apply retrogradation, then inversion and a transposition $t'=n-t=5$ to B, we return to the series $A$.

Let now $A$ be an arbitrary series and $B$ its RI with a transposition of t tone-fractions. Our observation tells us that applying IR with a transposition $n-t$ to $B$ yields $A$. Now let $\sigma$ be the PP obtained by applying RI with a transposition $t$ to $A$; we know that $\sigma$ maps $A$ into $B$. Now, since $\sigma$ maps $A$ into $B$, we know that the permutation that maps $B$ into $A$ is the inverse permutation of $\sigma$. Furthermore, by our observation, the permutation that maps $B$ into $A$ is the PP obtained by applying IR with the transposition $n-t$ to $B$. 
To sum up, the PP obtained from $A$ using RI is the inverse of the one from $B$ using IR: the proliferations will be the same except in that they will appear in inverse order, as we can see in Figure~\ref{fig: RIvsIR}, using the same $A$ as before:
\begin{figure}[htb]
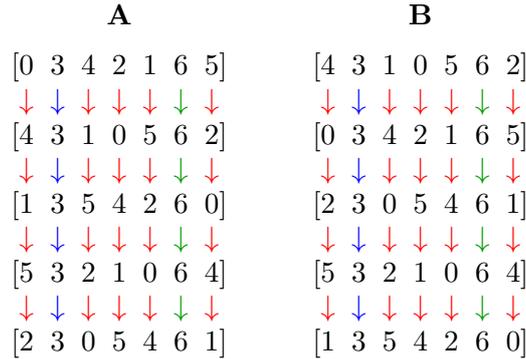

	
    \[
    \begin{array}{c@{\qquad}c}
    \textbf{A} & \textbf{B} \\[6pt]
    
    \begin{array}{c}
    [ 0\;\;3\;\;4\;\;2\;\;1\;\;6\;\;5 ] \\
    \textcolor{red}{\downarrow}\;\;\textcolor{blue}{\downarrow}\;\;\textcolor{red}{\downarrow}\;\;\textcolor{red}{\downarrow}\;\;\textcolor{red}{\downarrow}\;\;\textcolor{smoothgreen}{\downarrow}\;\;\textcolor{red}{\downarrow} \\
    \relax [ 4\;\;3\;\;1\;\;0\;\;5\;\;6\;\;2 ] \\
    \textcolor{red}{\downarrow}\;\;\textcolor{blue}{\downarrow}\;\;\textcolor{red}{\downarrow}\;\;\textcolor{red}{\downarrow}\;\;\textcolor{red}{\downarrow}\;\;\textcolor{smoothgreen}{\downarrow}\;\;\textcolor{red}{\downarrow} \\
    \relax [ 1\;\;3\;\;5\;\;4\;\;2\;\;6\;\;0 ] \\
    \textcolor{red}{\downarrow}\;\;\textcolor{blue}{\downarrow}\;\;\textcolor{red}{\downarrow}\;\;\textcolor{red}{\downarrow}\;\;\textcolor{red}{\downarrow}\;\;\textcolor{smoothgreen}{\downarrow}\;\;\textcolor{red}{\downarrow} \\
    \relax [ 5\;\;3\;\;2\;\;1\;\;0\;\;6\;\;4 ] \\
    \textcolor{red}{\downarrow}\;\;\textcolor{blue}{\downarrow}\;\;\textcolor{red}{\downarrow}\;\;\textcolor{red}{\downarrow}\;\;\textcolor{red}{\downarrow}\;\;\textcolor{smoothgreen}{\downarrow}\;\;\textcolor{red}{\downarrow} \\
    \relax [ 2\;\;3\;\;0\;\;5\;\;4\;\;6\;\;1 ] \\
    \end{array}
    &
    \begin{array}{c}
    [ 4\;\;3\;\;1\;\;0\;\;5\;\;6\;\;2 ] \\
    \textcolor{red}{\downarrow}\;\;\textcolor{blue}{\downarrow}\;\;\textcolor{red}{\downarrow}\;\;\textcolor{red}{\downarrow}\;\;\textcolor{red}{\downarrow}\;\;\textcolor{smoothgreen}{\downarrow}\;\;\textcolor{red}{\downarrow} \\
    \relax [ 0\;\;3\;\;4\;\;2\;\;1\;\;6\;\;5 ] \\
    \textcolor{red}{\downarrow}\;\;\textcolor{blue}{\downarrow}\;\;\textcolor{red}{\downarrow}\;\;\textcolor{red}{\downarrow}\;\;\textcolor{red}{\downarrow}\;\;\textcolor{smoothgreen}{\downarrow}\;\;\textcolor{red}{\downarrow} \\
    \relax [ 2\;\;3\;\;0\;\;5\;\;4\;\;6\;\;1 ] \\
    \textcolor{red}{\downarrow}\;\;\textcolor{blue}{\downarrow}\;\;\textcolor{red}{\downarrow}\;\;\textcolor{red}{\downarrow}\;\;\textcolor{red}{\downarrow}\;\;\textcolor{smoothgreen}{\downarrow}\;\;\textcolor{red}{\downarrow} \\
    \relax [ 5\;\;3\;\;2\;\;1\;\;0\;\;6\;\;4 ] \\
    \textcolor{red}{\downarrow}\;\;\textcolor{blue}{\downarrow}\;\;\textcolor{red}{\downarrow}\;\;\textcolor{red}{\downarrow}\;\;\textcolor{red}{\downarrow}\;\;\textcolor{smoothgreen}{\downarrow}\;\;\textcolor{red}{\downarrow} \\
    \relax [ 1\;\;3\;\;5\;\;4\;\;2\;\;6\;\;0 ] \\
    \end{array}
    
    \end{array}
    \]
    
	\caption[Example showing that the reasoning done for inversion and retrograding is equally valid for retrograding and then inversion.]{Comparison of the proliferations of A using RI with $t=2$ and B using IR with $t=5$}
	\label{fig: RIvsIR}
	
\end{figure}

In mathematical terms, we say that there is a bijection which preserves the proliferating structure between the PPs that use RI with $n$ notes and a transposition $t$, and the PPs that use IR with $n$ notes and the corresponding transposition $n-t$. This is why the structural properties are the same for RI and IR, so it will make no difference whether we use one or the other. From this point on, we will use the first option, that is, first inversion and then retrogradation. This will be convenient for notation, because this way we will only invert series that start from 0. With that concern out of the way, we will begin to examine how this kind of PPs work.

As we described in Section~\ref{sec: Prolif I}, the process of inversion+transposition creates pairs of elements which are related to each other: when we apply the process to one of them, we get the other one. We will call one number in this pair the inverse of the other one, for simplicity, even if this transformation is not only inversion; in reality, the existence of those pairs makes this transformation virtually an inversion. That is true except for the cases in which $x=-x+t$ modulo $n$, which will remain the same after applying the transformation. We will say those are their own inverses, or self-inverses. Those cases did not really matter for Section~\ref{sec: Prolif I}, but we will have to solve the equation for this one. We get  $2x=t$,  an equation which, using basic modular arithmetic, will have exactly one solution if $n$ is odd, in particular, $x=0$ if $t=0$, and a different number for every other transposition (we do not distinguish those as they will have the same properties). If n is even then it will have either two or no solutions according to whether t is even or odd, respectively. Those are exactly the cases that we will consider, each having slightly different results.

Retrogradation involves a completely new kind of behaviour: again, some form of pairs are involved, but now for the positions in the series rather than the numbers itself. The positions symmetric with respect to the center will be interchanged by the retrogradation. This means that, now, when we apply the PP to a given number $x$, we will not have the result of applying the transformation inversion+transposition to it, but rather we will get the result of applying that transformation to the note which is located as many spaces from the rightmost note of the series as the original is to the leftmost, which we will call its mirror.
Each position will be associated with its mirror, except for the case where n is odd, in which we will have one position, the center, whose mirror is itself. All other pairs of positions have the same behaviour; the only aspect that affects the proliferation is the number of pairs and if there is a center, so in that respect, they are completely interchangeable to us. Of course, there being a center position will change the behaviour of the proliferation, but we already had to distinguish between even or odd, so no new cases have to be considered. For simplicity, if we have a series of an odd number of notes and the PP that it generates, we will say that the cycle of the PP that contains the note of the series in this central position is the central cycle of the PP. 

If we focus on the first column from Figure~\ref{fig: RIvsIR}, we see that inversion+transposition transforms $A=[0, 3, 4, 2, 1, 6, 5]$ into $[2, 6, 5, 0, 1, 3, 4]$, so 0 is the inverse of 2, 3 is the inverse of 6, 4 is the inverse of 5 and 1 is self-inverse. With respect to retrogradation, 2 is its own mirror, 0 is the mirror of 5 and so on. We can factorise the permutation to obtain $\textcolor{red}{(}0 \; 4 \; 1 \; 5 \; 2\textcolor{red}{)}\textcolor{blue}{(}3\textcolor{blue}{)}\textcolor{smoothgreen}{(}6\textcolor{smoothgreen}{)}$, and, since the central note is 2, our central cycle is $\textcolor{red}{(}0 \; 4 \; 1 \; 5 \; 2\textcolor{red}{)}$.

\subsubsection{When $n$ is even and $t$ is odd}
\label{subsec: RI n even t odd}

Now we will explain in detail the case where n is even and t is odd, and all the others will have very similar reasonings. In this case, our equation $x=-x+t$  will have no solutions, and so, we will have exactly $n/2$ pairs of numbers which are mutual inverses. And also, we have no central position, because of n being even. 
The key observation, which will solve this case and help a lot for the others, is the following: the PP that we obtain after applying IR and an even number of transpositions to a series with an even number of notes cannot have two inverses in the same cycle. Furthermore, if the PP has some cycle $\textcolor{red}{(}a\; b\; c…\textcolor{red}{)}$ in its factorization, then $\textcolor{blue}{(}…c’\; b’\; a’\textcolor{blue}{)}$ also belongs to the PP. From here, we denote the inverse of any note $x$ by $x’$.

\begin{figure}[htb]  
	\centerline{
		\resizebox*{15cm}{!}{\includegraphics{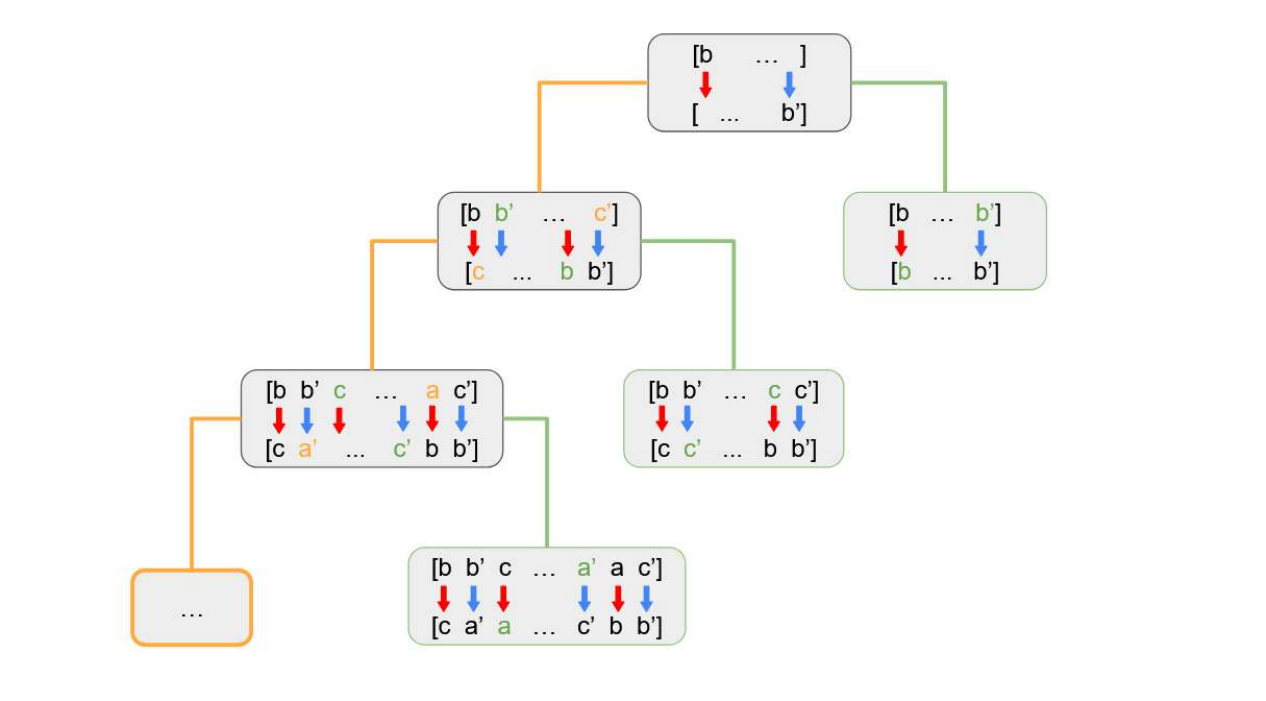}}\hspace{0pt}}
	\caption[Flowchart that shows how an arbitrary cycle is formed when n is even and t odd depending on the position of the notes in the series.]{Inductive process of the construction of a cycle when $n$ is even and $t$ odd}
	\label{fig: Closing RI cycle t odd}
	
\end{figure}

In order to show this, let us suppose that we have a note called $b$ in some position, and we will observe the cycle of the PP in which $b$ is involved. All the pairs of positions behave the same with respect to the PP, so, without loss of generality, we can assume that $b$ is in the first position. If $b'$ is in the opposite position, then $b$ will be sent by the PP to itself again, and the cycle will close. In any other case, $b’$ will be in some other position. Again, without loss of generality, we can assume it is in the second position. And $b$ will be sent to some other number $c$, meaning that $c’$ is in the opposite position of $b$. Again we could have $c$ be in the opposite position of $b’$, which would close the cycles $\textcolor{red}{(}b \;c\textcolor{red}{)}$ and $\textcolor{blue}{(}c’\; b’\textcolor{blue}{)}$. Or we could have another note a in the opposite position of $b’$, continuing the cycles as $\textcolor{red}{(}a\; b \;c…\textcolor{red}{)}$ and $\textcolor{blue}{(}c’\; b’\; a’…\textcolor{blue}{)}$. We can again suppose that $c$ is in the third position, and we will distinguish whether $a’$ is in its opposite position or not. If it is, the cycles close again, and if not, then another note $d$ is the image of $c$ through the PP… This process continues until the cycles close and, as we can see, neither $b$ nor any other note is ever in the same cycle as its inverse, which proves our observation inductively. This process is illustrated in Figure~\ref{fig: Closing RI cycle t odd}.

Now let us express this observation in terms of the obtainable PPs for this method. We know that, for any cycle in a PP, there has to be another cycle of the same length which is composed of the inverses of the notes in the first cycle. This means that all the permutations that could potentially be PPs have a structure which is the same as a permutation of $n/2$ notes but repeating every cycle 2 times.

At this point, we can notice that, in fact, every permutation of these characteristics is a PP. And we can construct each of them by just following the inductive process described above, as we can close the cycles any moment we want. For example, if we fix $n=24$ (even) and some $t$ odd, a permutation with the structure $[\{1\},\{1\},\{1\},\{1\},\{3\},\{3\},\{7\},\{7\}]$ is a valid PP since it repeats every cycle of the permutation $[\{1\},\{1\},\{3\},\{7\}]$ twice. A permutation given by the structure $[\{2\},\{2\},\{2\},\{9\},\{9\}]$ will not be possible, since it has an odd number of cycles of length 2.

\subsubsection{When $n$ is even and $t$ is even}
\label{subsec: RI n even t even}

Now if $n$ is even and $t$ is also even, our equation will have two solutions, so we will have two notes which are its own inverses. And there is still no center position in the series. Let $x$ and $y$ be the notes which are their own inverses, we could say that $x’=x$ and $y’=y$. 

\begin{figure}[htb]  
	\centerline{
		\resizebox*{15cm}{!}{\includegraphics{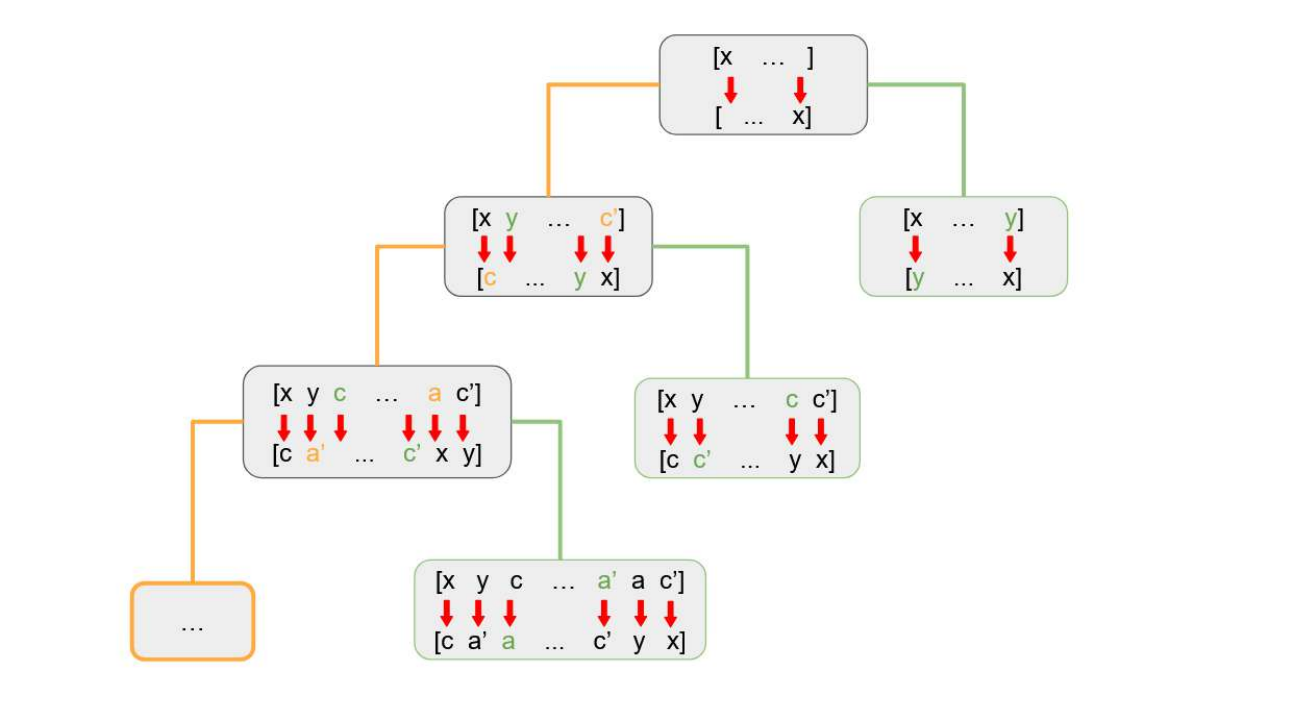}}\hspace{0pt}}
	\caption[Flowchart that shows how the cycle containing both self-inverses is formed when n and t are even.]{Inductive process of the construction of a cycle when $n$ is even and $t$ is also even}
	\label{fig: Closing RI cycle t even}
	
\end{figure}

Following the previous reasoning, but changing the role of b and b' for x and y, respectively, we start with x in the first position. Now, if y is its mirror, we will have the cycle $\textcolor{red}{(}x\;y\textcolor{red}{)}$, instead of $x$ and $y$ being sent to themselves. Otherwise, we assume that y is in the second position, so we could close the cycle taking the form $\textcolor{red}{(}x\;c\;y\;c'\textcolor{red}{)}$ instead of $\textcolor{red}{(}x\;c\textcolor{red}{)}$ and $\textcolor{blue}{(}c'\;y\textcolor{blue}{)}$. In general, our PP will have the cycle $\textcolor{red}{(}c’\; x\; c\; …\; e\; a\; y\; a’\; …\; f\textcolor{red}{)}$ instead of the previous $\textcolor{red}{(}a\; x \;c…e\textcolor{red}{)}$ and $\textcolor{blue}{(}c’\; y\; a’… f\textcolor{blue}{)}$. Figure~\ref{fig: Closing RI cycle t even} presents a diagram which shows how the cycle which contains the two self-inverses is constructed. 

In this case, any PP will be like a permutation of the previous case, but taking two cycles of the same length and combining them to get one of double the size. In other words, all the possible PPs for this case are obtained by choosing an even number $k \leq n$, which will be the length of one of its cycles (the cycle which contains the two self-inverses), and then any permutation of $(n-k)/2$ notes repeating every cycle of this permutation twice. One example of a PP for this case is given by the structure $[\{2\},\{2\},\{8\}]$ (notice that this is the permutation obtained in Figure~\ref{fig: barraque construction}). Here $k=8$, and we obtain our permutation repeating the cycles of $[\{2\}]$ twice and joining the result with \{8\}. The permutation $[\{6\},\{6\}]$ is not obtainable in this case.

\subsubsection{When $n$ is odd and $t$ is zero}
\label{subsec: RI n odd t zero}

\begin{figure}[htb]  
	\centerline{
		\resizebox*{15cm}{!}{\includegraphics{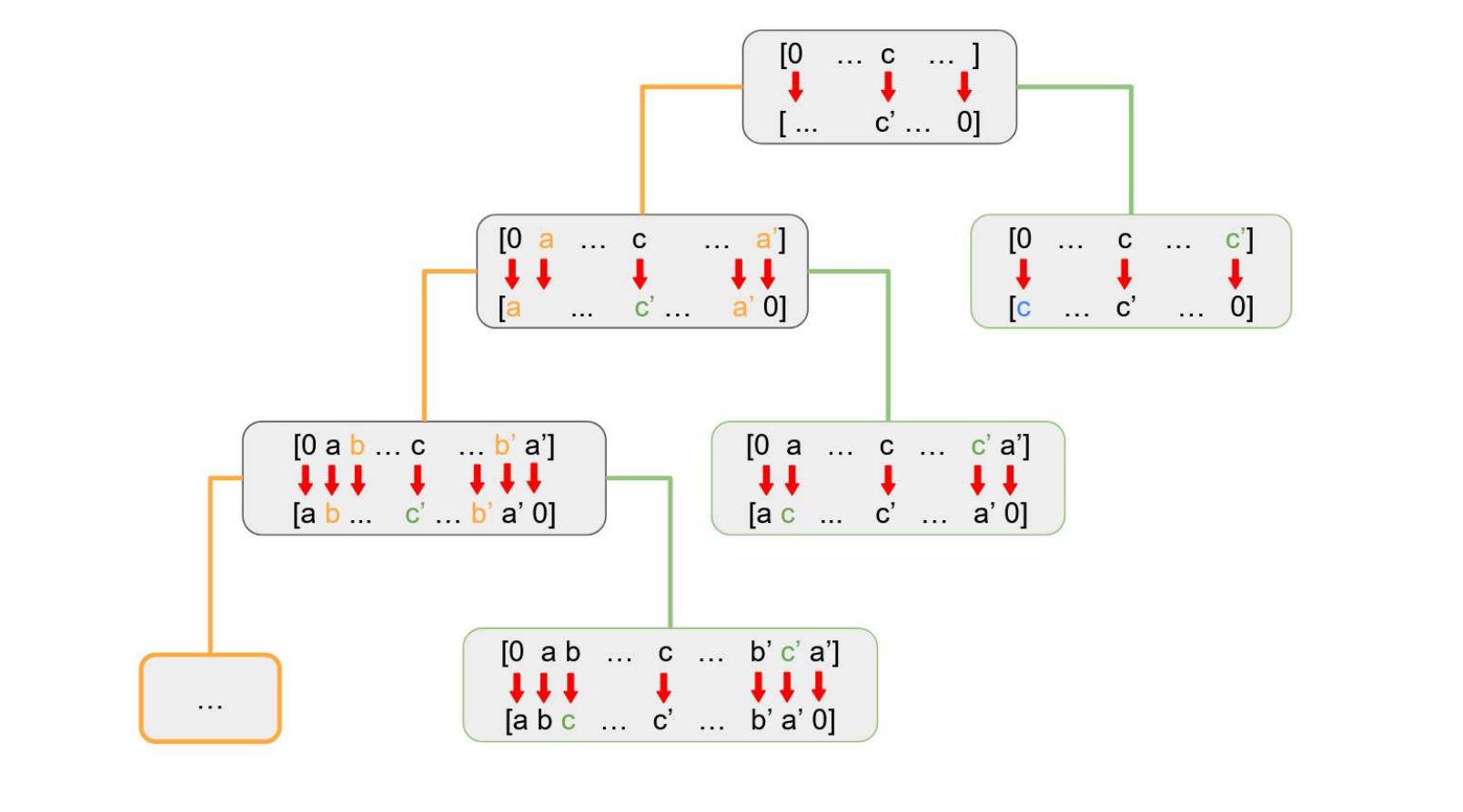}}\hspace{0pt}}
	\caption[Flowchart that shows how the cycle containing the only self-inverse is formed when n is odd and t is zero, using the center note to join both halves of the cycle.]{Inductive process of the construction of a cycle when $n$ is odd and $t=0$}
	\label{fig: Closing RI cycle n odd}
	
\end{figure}

In this case, there is only one number that is its own inverse, specifically the number zero, as was explained in the beginning of the section. Also, there is a central position in the series, meaning that the note in this place will be sent by the PP to its inverse
Applying the same reasoning from the previous section, the cycle that contains zero, the only self inverse, will have to look something like this: $\textcolor{red}{(}...\; b’\; a’\; 0\; a\; b\; …\textcolor{red}{)}$. In that case, this pattern continued until we hit the other self-inverse, closing the cycle this way: $\textcolor{red}{(}c’\; y\; c\; …\; b’\; a’\; 0\; a\; b\; …\textcolor{red}{)}$. Now this cannot happen, since we do not have another self-inverse. But now the central position plays a role, since it is the only other change in our structure and, therefore, the only way in which the cycle can close. We continue in the manner displayed in Figure~\ref{fig: Closing RI cycle n odd} until we reach the central note (call it $c$), which now will not be mapped into some other arbitrary note, but rather its own inverse, $c’$. This will close the cycle in this form: $\textcolor{red}{(}c’\; …\; b’\; a’\; 0\; a\; b\; …\; c\textcolor{red}{)}$.

Thus, the cycle containing the note zero has now an odd length greater than one, since 0 is in the first position of our series so it cannot be mapped to itself (except when $n=1$). All other cycles behave as they did in previous cases. That is to say, the possible PPs in this case are those which have a cycle of length $k$, for some odd number $k$ greater than 2,  and then any permutation of $(n-k)/2$ notes repeating every cycle of this permutation twice. For $n=9$ (and $t=0$), $[\{2\},\{2\},\{5\}]$ is a PP, using $k=5$, but $[\{4\},\{5\}]$ is not.

\subsubsection{When $n$ is odd and $t$ is not zero}
\label{subsec: RI n odd t nonzero}

This last case has once more only one number that is its own inverse, but now, any one of the notes in our series (except zero) can hold this property depending on the transposition as opposed to before, where it could only be zero. We also have a central position because of $n$ being odd, so the only difference with the last case is that the self-inverse can be in any position.
Any other configuration of the previous case can be achieved by putting the self-inverse in a position that is not central, and we get a new case: now the cycle containing the self inverse can be of length 1, if it is in the center.
Thus, the cycle containing the self-inverse has an odd length, and all other cycles behave as they did in Section~\ref{subsec: RI n odd t zero}. That is to say, the possible permutations obtainable in this case are those which have a cycle of length $k$, for some odd number $k$ (now it can be 1), and then any permutation of $(n-k)/2$ notes repeating every cycle of this permutation twice. For $n=15$ (and $t\neq 0$), $[\{1\},\{7\},\{7\}]$ is a PP, using $k=1$, but $[\{2\},\{2\},\{3\},\{8\}]$ is not.

\subsection{Study of proliferations using retrogradation (R)}

The final case is the most complex one because, in some way, the variety of new possibilities according to the possible values of $n$ and $t$ generalizes the structures that we had in the other sections. In this case we will have to deal with retrogradation and transposition, but not inversion. Ironically, inversion makes the proliferating series simpler, as it forces pairs of inverses, which make the analysis easier. As in Section~\ref{sec: Prolif RI}, we will have to treat transposition and retrogradation as separated processes: first, the transposition will send each note $x$ to $x+t$, as in Section~\ref{sec: Prolif P}, and then we will make the retrogradation, so that the PP sends one note $x$ to the transposed of its mirror.

As we discussed in Section~\ref{sec: Prolif P}, the transposition will create $gcd(n,t)$ cycles of length $n/gcd(n,t)$. Note that those cycles are not the cycles in our PP, as we have not applied retrogradation yet, so we will have to use the term “cycle” with those two different meanings: those of the transposition and those which form the PP itself. Again, as in Section~\ref{sec: Prolif P}, the only way in which t affects our PP is the quantity and length of the transposition cycles, and those depend only on the greatest common divisor. That means that, given a number of notes $n$, different values of t with the same greatest common divisor with $n$ will have exactly the same results: the same PPs will be possible and the same quantity of series will lead to a given PP. The only factor that could possibly change is the specific series which generate the PP. To illustrate this, if we choose $n=12$, this reasoning justifies that we will obtain the same results for $t=2$ and $t=10$, as $gcd(2,12)=gcd(10,12)=2$. 

In other words, given a value of $n$, what really determines how the proliferation works with respect to $t$ is the number of transposition cycles and their length. This ultimately depends, as we said, on $gcd(t,n)$, but thinking about the transpositions just in terms of number and length of transposition cycles will help understand what we call the generalized transpositions, the tool that we will use for the general case of the proliferating series through retrogradation. Like in the last example, if we choose $n=12$ and $t=2$, our transposition cycles will be $(0\; 2\; 4\; 6\; 8\; 10)$ and $(1\; 3\; 5\; 7\; 9\; 11)$. Alternatively, for $t=10$, we get $(0\; 10\; 8\; 6\; 4\; 2)$ and $(1\; 11\; 9\; 7\; 5\; 3)$. Both of them generate 2 ($gcd(n,t)$) cycles of length 6 ($n/gcd(n,t)$).

We will divide the explanation of the proliferating series through retrogradation into four sections: first, two particular cases which behave quite simply compared to the rest, for $t=0$ and $t=n/2$. A third section for the case of $t$ and $n$ being coprime, that is to say, $gcd(n,t)=1$, a very illustrative case for the ideas of the fourth section, and finally a section for the general case, which encompasses every possibility, including the previous cases.

\subsubsection{When $t$ is zero}
\label{subsec: R t zero}

This particular case is of little interest because it produces no variety whatsoever, the same as in Section~\ref{sec: Prolif I}. The retrograding process sends the first note to the last, the second to the second last, etc. In short, the PP relates each note to the one in its mirror position, resulting in cycles which are always of length 2 (except in the case of the note in the middle position when $n$ is odd, which produces a cycle of length 1). Again, the global PP will have order at most 2, and will not generate new series.

\subsubsection{When $n$ is even and $t=n/2$}
\label{subsec: R t is n/2}

Given an even value of $n$, the case of $t=n/2$ is remarkable for its relation with the previous reasonings. If we analyze the transposition, we always have $gcd(n,t)=n/2$, so the transposition is made up of $n/2$ cycles of length 2. Those cycles of length 2 are pairs of notes which are sent to each other once we apply the transposition. And then we have to apply the retrogradation to get the PP. Recall from Section~\ref{sec: Prolif RI}: we have a transformation, inversion+transposition, which generates pairs of inverses and up to 2 fixed points, depending on the case. Then we have to apply retrogradation to get the PP. Our transformation is now just transposition, but it indeed consists of pairs of “inverses”, with no fixed points, which has the same structure as RI with $n$ even and $t$ odd, because this transposition behaves the same as a regular inversion in terms of cycles. This concludes that the possible PPs obtainable for this case, for a fixed even $n$ and $t=n/2$ with retrogradation, are exactly the same as in the case of retrograde inversion for the same even $n$ and any $t$ odd. Not only that: any of those PPs that we choose is achieved through retrograde inversion (with $t$ odd) with the same number of series as through just retrogradation (with $t=n/2$); it is just a matter of labeling, but the structural properties are the same.

\subsubsection{When $gcd(n, t) = 1$}
\label{subsec: R gcd 1}

Let us suppose that $gcd(n, t)=1$. For simplicity's sake, the proof will be written as if $t=1$ even though the same reasoning applies for any transposition that meets this condition. It is also worth noting that the general case described in Section~\ref*{subsec: R general} would suffice in order to justify the way in which the series here behave. However, this previous case study has some concrete results that cannot be generalized, and it will also help understand the ideas of the general case, so we have reason to first examine the case where $gcd(n,t)=1$.

The way in which we study the series shall be through a recursive process which will help us determine the possible combinations of cycles that can result from any series of a given number of notes: the idea is, supposing that we know all the possible PPs for a given number of notes $n$, with some transformation (adding a note to the series), to extend those to all the possible PPs of $n+1$ notes. This way, we start from a base case that we already know for its simplicity (in this case $n=2$), and extend it recursively to achieve every possible value of $n$.  

For $n=2$, the only possible structure is clear (recall $t=1$):

\begin{tabular}{l c}
Original & → [0, 1] \\
 &  \;\;\; $\textcolor{red}{\downarrow} \;\; \textcolor{blue}{\downarrow}$\\
Retrog + transp & → [0, 1] \\
\end{tabular}

Now suppose that for a certain number $n$ of notes, we know all the possible structures that the PPs can have, and let's use this to figure out the structures that can be obtained when dealing with $n+1$ notes. We will distinguish two separate cases: the first, where $n$ is an even number, and therefore $n+1$ is an odd number, and the second, where $n$ is an odd number and $n+1$ is an even number.

Let us start with $n$ being an even number and consider an arbitrary series with $n$ notes. To make sure that our inductive process captures every possible PP of $n+1$ notes, we have to check that every PP of $n+1$ notes can be obtained from a series of $n$ notes applying some transformation to extend it to $n+1$ notes. In this case, given a series of $n+1$ notes, we can just eliminate the note in the center, resulting in a series of $n$ notes such that, when applying transposition, we just skip that note, and all the retrograding relations will maintain, yielding a PP with similar structure. Thus, our inductive step from $n$ to $n+1$ will consist of adding a new note in the center position to a series of $n$ notes. The name of this new note itself is of no particular relevance, but given that we already have notes from 0 to $n-1$ in the series of $n$ notes, we shall call this new note $n$. Before, our transposition sent $n-1$ to 0, since we established $t=1$. Now we are doing arithmetic modulo $n+1$ instead of $n$, so $n-1$ will be sent to $n$, and $n$ to 0. Every other relation will be maintained.

With respect to the PP, our new note $n$ will be sent to 0, since it is in the central position, and the note in the mirror position of the one numbered $n-1$, which in the n-note series was paired with 0, now gets sent to our newly added note $n$. The result of this for our n+1-note series is a cycle structure which is identical to that of the n-note series that it derives from, except for the cycle that contains the note 0 (or any other cycle, since we can relabel the notes), in which one more note will be added:

\[ \left[ \begin{array}{cccc}
	n-1, ... & ..., b \\
    \textcolor{red}{\downarrow} \;\;\;\;&\;\;\;\; \textcolor{blue}{\downarrow}\\
	b+1, ... & ..., 0 
\end{array} \right]
\longrightarrow
\left[ \begin{array}{ccc}
	n-1,  …,&  n,& … , b \\
    \textcolor{orange}{\downarrow}\;\;\;\;\;\;\;\;&\;\textcolor{smoothgreen}{\downarrow}\;\;&\;\;\;\;\;\;\textcolor{smoothgreen}{\downarrow}\\
	b+1, …,    &0,& …,   n
\end{array} \right]
\]

Now let $n$ be an odd number. Once again, we have to obtain an inductive process that assures us that every PP of $n+1$ notes can also be obtained from some inductive extension of a series of $n$ notes. In this case, because of $n+1$ being even, we can eliminate any note to get a series of $n$ notes in which every relation caused by retrogradation is maintained except for the mirror of the eliminated note, which now will have no mirror note and will not be affected by retrogradation. Then, our inductive step will take every series of $n$ notes, and add a new one, again called $n$, in a symmetrical position to the one that was previously the note in the center. For the following reasoning, we shall call the note that was in the center in the series of $n$ notes $a$.

\[[x_1, x_2, …, a,  …,  x_n]   \longrightarrow  [x_1, x_2 …,  a, n,  …,   x_n]\]

At this point, in a way that is similar to what was obtained in the previous case, the mirror of $n-1$ goes to $n$ when applying the permutation, while $n$ will go to $a+1$, as we can see below.  

So, how does the new note $n$ modify the cycle structure when going from $n$ to $n+1$? In this case, two possible outcomes can be achieved: 

\begin{enumerate}
	\item[(i)] Two cycles get joined together and an additional note is added to the result.
	\item[(ii)] A cycle is broken into two and a note is added to one of the resulting cycles. 

\end{enumerate}

This depends on whether or not the notes that we named a and 0 are in the same cycle. If they are, the cycle which contains them both is broken. If they are not, the cycles containing one note and the other join together.

\[ \left[ \begin{array}{ccc}
	n-1,  …, &  a,&  …, b \\
    \textcolor{red}{\downarrow}\;\;\;\;\;\;\;&\textcolor{blue}{\downarrow}\;&\;\;\;\;\;\;\textcolor{darkpink}{\downarrow}\\
	b+1, …,& a+1,& …, 0 
\end{array} \right]
\longrightarrow
\left[ \begin{array}{cccc}
	n-1,  …, & a, &   n,  …, &  b \\
    \textcolor{orange}{\downarrow}\;\;\;\;\;\;&\textcolor{smoothgreen}{\downarrow}&\textcolor{purple}{\downarrow}\;\;\;\;\;\;\;&\textcolor{purple}{\downarrow}\\
	b+1, …,  & 0, &   a+1, …, & n
\end{array} \right]
\]

In the left series displayed above we have a cycle  of the form $\textcolor{darkpink}{(}b, 0, …B_1…\textcolor{darkpink}{)}$ and another one    $\textcolor{blue}{(}a, a+1, …B_2…\textcolor{blue}{)}$, where $B_1$ and $B_2$ represent the unknown part of the cycles. 
It is possible that the permutation under consideration has $a$ and 0 in the same cycle. In this case we obtain the cycle $\textcolor{darkpink}{(}b\; 0\; …B_1…, a, a+1, …B_2…\textcolor{blue}{)}$, rather than the two cycles mentioned before.  The series on the right  in the cycles represented shows that $a$ now gets sent to 0 instead of $a+1$, $b$ is sent to our new note $n$ instead of 0, and $n$ is sent to $a+1$. Every other note in the PP will behave the same, so the blocks $B_1$,$B_2$ after 0 and $a+1$ will remain unaltered.

Thus, in the first case, where the cycles which have 0 and $a$ are different, those will join in the series in the right to form the cycle $\textcolor{purple}{(}b, n, a+1, …B_2…, a, 0, …B_1…\textcolor{smoothgreen}{)}$. 
In the second case, the cycle which contains both 0 and $a$ will separate into $\textcolor{purple}{(}b, n, a+1, …B_2…\textcolor{purple}{)}$ and $\textcolor{smoothgreen}{(}a, 0, …B_1…\textcolor{smoothgreen}{)}$.

To sum up, the inductive process works like this: we start from $n=2$ knowing that the only possible PP is the identity permutation. Now let $n$ be a natural number and suppose that we know all the possible PPs of $n$ notes. Then, our inductive step brings us all the possible PPs of $n+1$ notes, which are described as follows:

\begin{enumerate}
	\item[(i)] If $n$ is even, then we take every PP of $n$ notes. For each cycle of each PP, we have another PP of $n+1$ notes extending that cycle by 1 note.
	\item[(ii)] If $n$ is odd, once again we take every PP with $n$ notes and the center cycle. We get a new PP of $n+1$ notes for every possible partition into two cycles of that cycle, adding an extra note in one of them. And also another new PP for each other cycle of that PP that we can combine with it, adding one extra note.

\end{enumerate}

With this in mind, we can use the inductive process to establish a pretty strong description that tells us every PP possible for each $n$ without need of the recursive process more than for proving the description true. For that, we have to make two observations:

On the one hand, let us observe the number of cycles that a PP of $n$ notes can have. For $n=2$, our only PP is the identity, which has two cycles of length one. The inductive step from $n$ even to $n+1$ odd cannot add new cycles to the structure, so $n=3$ will also have PP with two cycles. The step from $n$ even to $n+1$ odd may break up a cycle or combine two, that is, we will have PPs with one more or one less cycles, so the maximum number of cycles of PPs of four notes is three. Extending this reasoning for every $n$, we get that, if it is even, then the maximum number of cycles is $(n/2)+1$, and if it is odd, $(n+1)/2$. In a more compact way, we can say that the number of cycles is always less than or equal to $(n/2)+1$. We will call this \textbf{condition 1}.

On the other hand, let us observe the number of cycles of even length that a PP can have. Given a PP, let us call that number of cycles $k$. For $n=2$, as we said, the only PP has two cycles of length 1, so this PP has $k=0$, an even number. 
When going from $n$ even to $n+1$ odd, our inductive process consists of adding 1 length to some cycle of the PP. If that cycle had an even length, now it will have an odd length, so $k$ will decrease by 1. If, by contrary, it had an odd length it will now have an even length, increasing $k$ by one. Going from $n=2$ to $n=3$, our only PP had $k=0$, so, since a permutation cannot have $k=-1$, every PP of 3 notes will have $k=1$. Observe that decreasing or increasing $k$ by one always changes its parity.
In the other case, going from $n$ odd to $n+1$ even, we join two cycles or split one into two, and then add 1 to the length of one of the remaining cycles. When joining two cycles, we have three options: if both are even, then we combine them and add one note to the result, losing two cycles of even length to create a new one of odd length, decreasing $k$ by 2. If one is even and the other one is odd, we lose one even cycle but create a new one of even length, leaving $k$ invariant. If both are odd, then the combined cycle is also odd, so $k$ also remains invariant. When splitting a cycle into two, the results are reversed: $k$ can increase by 2 or remain the same. Observe that all of those operations, increasing or decreasing by 2 and leaving it the same, always conserve the parity of $k$. 

This tells us that the parity of $k$ is invariant for all the PPs of a fixed number of notes. That is to say, given a number of notes $n$, different PPs may have different values of $k$, but always sharing the parity. In fact, we know that for $n=2$, we have $k=0$, which is even, for the only PP. For $n=3$ all PPs will have $k=1$ (odd), for $n=4$ all will have $k$ also odd, for $n=5$ it will change again to $k$ even, and so on. We get the following result: every PP obtained from a proliferating series will have an even value of $k$ associated if $n=1,2$ modulo 4, and an odd value of $k$ if $n=0,3$ modulo 4. We will call this \textbf{condition 2}.

With that, we have two conditions that every PP must follow: the maximum number of cycles that it can have and the parity of its associated $k$. It turns out that any permutation that meets just those two conditions is indeed a PP, a remarkable result which we will prove next. We will in fact prove a slightly stronger result, adding the fact that, if $n$ is odd and $\sigma$ is a PP, any cycle of $\sigma$ that is not of length one can be the central cycle. To prove this by induction, we check that it holds for $n=2$ and prove that if it holds for all numbers up to $n$, then it holds for $n+1$.

There are just two permutations of 2 notes: the one with two cycles of length one and the cycle of order two. The first one meets conditions 1 and 2, but the second does not satisfy condition 2. As we also saw, the first one can be obtained through proliferation and the second can’t, so this shows that the result holds for the base case. Now suppose that $n \geq 2$ and the result is true up to $n$ and let $\sigma$ be a permutation of $n+1$ notes with the structure $[\{C_1\},\{C_2\},...\{C_m\}]$ that satisfies conditions 1 and 2. Recall that \{$C_i$\} represents a cycle of length $C_i$. We will see that there exists a PP of $n$ notes which we can extend to $\sigma$ with our inductive transformations, proving that there exists a series of $n+1$ notes that generates $\sigma$, that is, $\sigma$ is a PP. Note that $(n+1)/2$ and $(n/2)+1$, the maximum numbers of cycles that condition 1 allows, are strictly less than $n$ for $n \geq 3$, and since we are reasoning for $n \geq 2$, then $n+1 \geq 3$. Then, the number of cycles that $\sigma$ can have must be strictly less than $n+1$, meaning that there is some cycle that has length at least 2. Without loss of generality, suppose that the length of $\{C_1\}$ is strictly greater than 1. The inductive step will once again be separated when $n$ is even and when $n$ is odd:

If $n$ is even, consider a PP of $n$ notes $\tau$ given by the structure $[\{C_1-1\},\{C_2\},...,\{C_m\}]$. It satisfies condition 1 because $n/2+1=((n+1)+1)/2$. And it satisfies condition 2 because $\sigma$ satisfies it and decreasing $\{C_1\}$ by one changes the parity of $k$, which is what must happen in this case. Since $\tau$ satisfies the conditions, there must be by induction a series that produces $\tau$ with proliferation. And decreasing $C_1$ by 1 is the opposite of our inductive transformation from $n$ even to $n+1$ odd, so we can extend that series to a series of $n+1$ notes that produces precisely the permutation $\sigma$. The only thing left is checking that the central cycle can be any of the cycles of $\sigma$. But when extending a series in this way, the cycle which we increased is precisely the cycle in the center, so we could have done exactly the same process for any cycle of $\sigma$ of length at least 2 to obtain $\sigma$ from a series that has another cycle in the center. 

Now suppose that $n$ is odd. If $m>1$, that is, if $\sigma$ has more than one cycle, then consider the permutation of $n$ notes $\tau$ determined by $[\{C_1+C_2-1\},  \{C_3\},… ,\{C_m\}]$. The permutation $\tau$ has one less cycle than $\sigma$, so it satisfies condition 1, and for a reasoning similar to others exposed before, it also satisfies condition 2. Then, by induction, $\tau$ is a PP, and we can also suppose that  the cycle $\{C_1+C_2-1\}$ is the cycle in the center, since we can choose any cycle of length greater than 1. To extend this series to another one of $n+1$ notes, we will break the center cycle into $\{C_1\}$ and $\{C_2\}$, resulting in a series of $n+1$ notes that proliferates into $\sigma$, which is what we wanted. If $m = 1$ then $\sigma$ only has one cycle, so its structure is just $[\{C\}]$, where $C = n+1$. Since $n$ is odd and $n \geq 2$, we can assure that $n \geq 3$. For this case, consider $\tau$ given by the structure $[\{C_1\}, \{C_2\}]$, where $C_1 = n-1$ and $C_2 = 1$, which satisfies condition 1 since it only has two cycles and satisfies condition 2 because $n-1$ has the same parity as $n+1$, so $k=1$ for both $\sigma$ and $\tau$. Then, there again exists a series of $n$ notes which proliferates into $\tau$, which can have $\{C_1\}$ in the center since $n-1 \geq 2$. We extend this series to $n+1$ notes joining the cycles $\{C_1\}$ and $\{C_2\}$ to obtain a series that produces $\sigma$ in the proliferation process. 

We started from a permutation $\sigma$ of $n+1$ notes satisfying conditions 1 and 2 and, by induction, constructed a series of $n+1$ notes that proliferates into that permutation, meaning that $\sigma$ is a PP, which concludes the proof.

\subsubsection{General case: generalised transpositions}
\label{subsec: R general}

The process through which we obtain the possible structures for proliferating series of any number of notes $n$ and any transposition $t$ follows a thought process that is similar to that of Section~\ref{subsec: R gcd 1}.

In order to generalise the previous section, we need to define a concept which will be key in the study at hand: generalised transpositions.
When we transpose a series of $n$ notes, the cycles of the permutation which results from said transposition are all of the same length. More specifically, we get $gcd(n,t)$ cycles of length $n/gcd(n,t)$. However, there is one problem that arises: when we add or subtract notes from the series while doing the induction process, the structure of our transposition is broken and thus the reason for this new concept. 

A generalised transposition (GT) is, \textit{a priori}, any transformation (any permutation) which will be used in place of a regular transposition to then be combined with retrograding as one would normally do. Since we are going to add notes one by one during the induction process, we want to start from regular transpositions and add a note in a single cycle of the transposition, repeatedly, to get different kinds of possible transformations. Adding one note in just a cycle means that, in this cycle, there was some note $a$ transformed into some other note $b$, that now has an intermediate note $c$, so that $a$ is transformed into $c$ and $c$ into $b$. And this includes adding a new note which is not related to any of the previous notes by the transformation, so it gets transformed into itself. Therefore, not all the structures studied in the induction process correspond to actual proliferating series constructed using the processes of traditional serialism, but they are all necessary steps towards building series with larger numbers of notes. 

For the base case in this study we will use $n=1$, which has a single possibility, since there exists only one GT, which is the identity permutation. Once again, we will start by going from an even number of notes to the next number, which will be odd, and afterwards we will do the opposite, going from an odd number $n$ to an even number $n+1$.

Suppose that $n$ is an even number for which we know all the possible PP structures given any GT.  Let’s consider a series of $n+1$ notes in which the note that is found in the central position is named $n$ (remember that this is not a restricting requirement because we can relabel all the notes in such a way that this happens without modifying the structure of the PPs). Associated with that series, we will have a certain GT involving all the notes.

Now, the structures that this n+1-note series might have are unknown to us, so we must relate it to a series with n notes, for which we do know the possible structures, and then work our way back to our n+1-note series. 
Going one step back into the n-note series is straightforward: simply remove the note $n$ which is in the central position. The n-note series we obtain has a GT associated which comes from removing a note in the GT that we had for the n+1-note series. So now we have an n-note series and a GT associated with it. Because of the induction, the cycle structure of this series after we apply the GT and then retrograde it is known. 

To figure out the structure of the n+1-note series, let us add the note $n$ again. If the new note is not a fixed point of the GT, the reasoning is similar to the case where $gcd(n,t)=1$, we need to add this note into one of the cycles of the PP that was given to us in the n-note series, but we now have one small restriction to take into consideration. The note $n$ can only be added to a cycle of the PP which contains a note that is in the same cycle as $n$ in the GT. This is consistent with the case study from Section~\ref{subsec: R gcd 1} because, with $n$ and $t$ being relatively prime, the transposition that is considered has a single cycle which contains all the notes. The restriction happens because, since $n$ is in the middle position of the series, the note to which it maps in the PP is the same note to which it goes in the GT. Therefore, all the series that stem from our chosen n+1-note series with the transposition will have in common that $n$ is in the same cycle as one of the notes with which it shares a cycle in the GT.

The last option is that $n$ is a fixed point of the GT. In that case, our PP will also have $n$ as a fixed point, and all the rest remains the same, as we can see here:

\[ \left[ \begin{matrix}
	 x &   … &    y  \\
    \textcolor{red}{\downarrow}&&\textcolor{blue}{\downarrow}\\
	 GT(y) & … & GT(x)
\end{matrix} \right]
\longrightarrow
\left[ \begin{matrix}
	 x &  … &   n&  … &  y  \\
    \textcolor{orange}{\downarrow}&&\textcolor{smoothgreen}{\downarrow}&&\textcolor{purple}{\downarrow}\\
	 GT(y) & … & n & … & GT(x) 
\end{matrix} \right]
\]

Suppose now that $n$ is odd and let us study what happens to an arbitrary series with $n+1$ notes and a certain GT of our choosing. As in Section~\ref{subsec: R gcd 1}, we will place $n$ right next to its mirror so that, when we eliminate $n$, we will not have to move around several pairs of notes, but instead keep the mirror of $n$ as the center note without a pair. 
We obtain a series of $n$ notes with a GT that is just like before, except for one cycle, which will be missing a note. 

When we put $n$ back in its place, doing a similar reasoning as in Section~\ref{subsec: R gcd 1}, three possible outcomes will be: 

\begin{enumerate}
	\item[(i)] The central cycle is divided into two and one of them gets added $n$ as an additional note.
	\item[(ii)] The central cycle and another cycle which contains a note that is in the same cycle as $n$ in the GT join together and $n$ is added to the resulting cycle.
	\item[(iii)] In the case that $n$ is a fixed point in the GT, the central cycle simply gets added one note : 
	
\end{enumerate}

	\[ \left[ \begin{matrix}
		x&  …&  a &  … &  y \\
        \textcolor{red}{\downarrow}&&\textcolor{blue}{\downarrow}&&\textcolor{darkpink}{\downarrow}\\
		GT(y) & … & GT(a) & … & GT(x)
    \end{matrix} \right]
	\longrightarrow
	\left[ \begin{matrix}
	x &  … &     a &  n &  … & y \\
    \textcolor{orange}{\downarrow}&&\textcolor{smoothgreen}{\downarrow}&\textcolor{smoothgreen}{\downarrow}&&\textcolor{purple}{\downarrow}\\
	GT(y) & … &   n &  GT(a) & … & GT(x)
	\end{matrix} \right]
	\]

As we can see, there are some modifications in the inductive process that we deduced in the last section. Now, our conditions 1 and 2 can be reformulated to something similar but also depending on the particular GT. But the further restrictions make it hard to find a closed mathematical description for these cases, so we only have an inductive process to construct them, which familiarizes us with these kinds of objects. This was also useful to make a program that calculates every possible PP for every $n$ and $t$, in Appendix A.

\section{Classification of proliferating series using RI}
\label{sec: Classification}

As we explained throughout Section~\ref{sec: Prolif n notes}, there are some structural properties that retrogradation and inversion+transposition induce. Those structural properties consisted in the fact of there being pairs of either inverses or symmetrical positions, and their exceptions (self-inverses and central position).  Then, we can deduce some transformation for series that do not affect the PP structure that they generate. Those consist in interchanging:

	\begin{enumerate}
	\item[(i)] Two different pairs of positions
	\item[(ii)] The two notes in a pair of positions
	\item[(iii)] Two different pairs of inverses
	\item[(iv)] The two notes in a pair of inverses
\end{enumerate}

We can also naturally extend those if we are not inverting for the cycles of length $n/gcd(n,t)$, using the same idea for cycles of length greater than 2. Those are all the operations that constitute the classification. In fact some of these operations can be reproduced using a combination of other operations, but that is not of interest for our analysis, and keeping all of them helps to better understand the intuition behind the ideas. The objective of this section is using those to classify the series with respect to the PP structure that they generate. The classification makes little sense for the cases of I and P, because the possible structures are determined just by $n$ and $t$, not the concrete series, but it will bring new information for RI and R.

Having explained the operations that leave the PPs invariant, the next step is grouping all the series that can be transformed into one another using those operations in what we will call an \textit{equivalence class}. If we fix a serial transformation for the proliferating process, and values of $n$ and $t$, the different equivalence classes represent the different PP structures for that transformation. This is, if two series are in the same class, then we can repeatedly apply our operations to any of them to obtain the other one, and the structure of the PPs obtained will be the same for both of them. What we would also want for this kind of classification is the reciprocal: two series being in different classes means that they will generate PPs with different structures. That is, given a PP structure, all of the series that proliferate into that PP will be in the same equivalence class. 

This reciprocal is not true, though, for the case of R. In the context of proliferating series using R, there are examples of values of $n$ and $t$ with which we can find different series proliferating into PPs of the same structure, but whose PPs have essential differences which make them be in different equivalence classes with this classification. One of those examples is shown in the next figure, using $n=8$ and $t=1$:

\[[0, \;1,\; 2,\; 3,\; 5,\; 4,\; 7,\; 6],\; [0,\; 1,\; 2, \;3, \;5, \;7, \;4 ,\;6]\]

Those two series proliferate into permutations with the structure of a single cycle of length 8, but cannot be obtained from each other by applying our transformations.

However it is true for the case of RI, so we can use this to classify those types of series with respect to the structure of PPs. This doesn’t invalidate this whole idea for R, but the classification would have to take into account parameters that are of little interest for proliferation, so we just have to consider that there are several classes for each PP structure in the case of R.

\begin{figure}[htb]
	\centerline{
		\resizebox*{18cm}{!}{\includegraphics{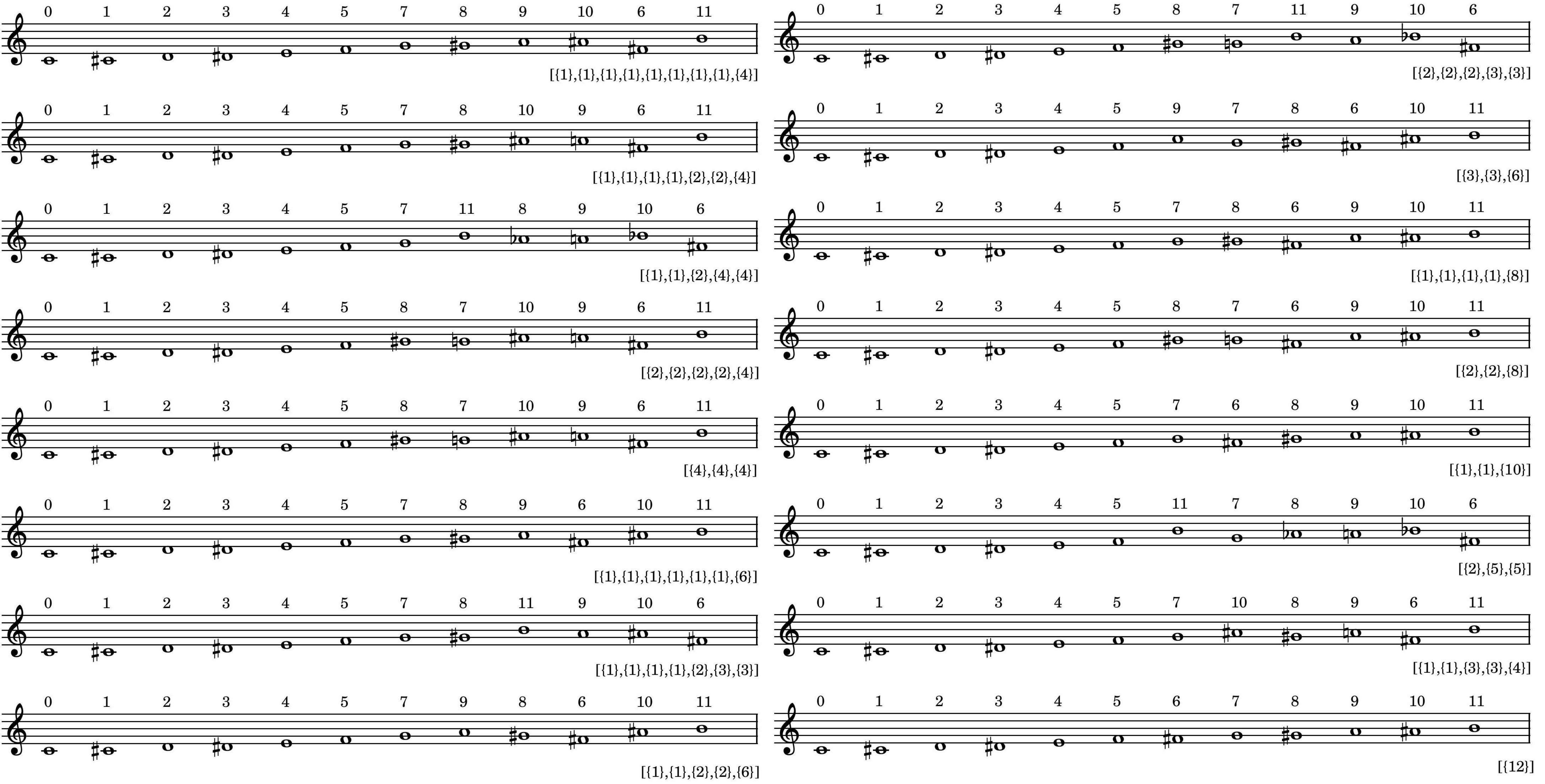}}\hspace{0pt}}
	\caption[List containing one element from each equivalence class. By applying the permitted transformations to each one of these elements we can obtain all other elements belonging to the same class.]{Base series for each obtainable cycle structure applying RI for $n$=12 and $t$=0}
	\label{fig:Classification}
	
\end{figure}

\begin{figure}[htb]
	\centerline{
		\resizebox*{18cm}{!}{\includegraphics{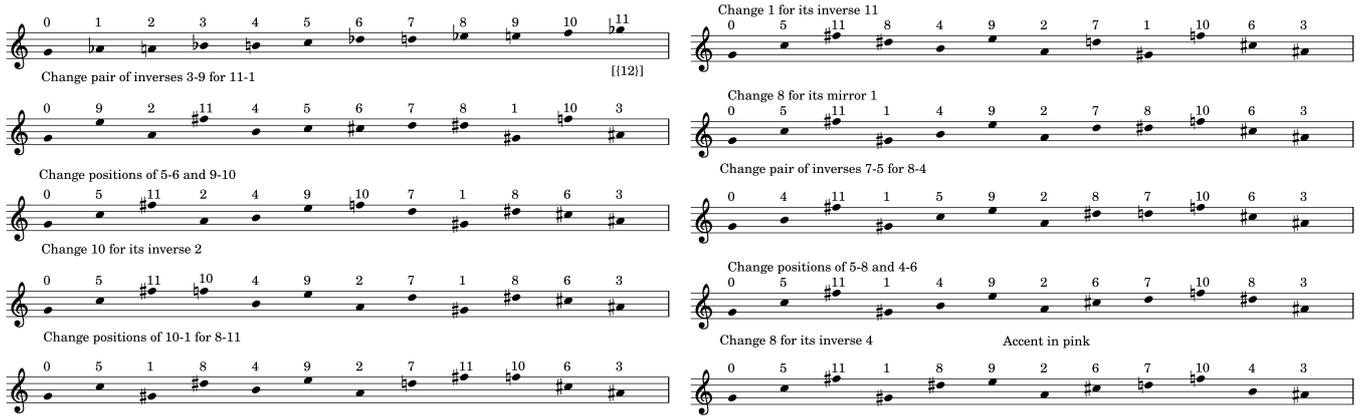}}\hspace{0pt}}
	\caption[Example of creation of a specific series with the desired properties starting from one of the base series in the classification.]{Use of the classification to create a proliferating series}
	\label{fig:Classification_use}
	
\end{figure}

Joining all this together, what this reasoning allows us to do is picking one series for each PP structure obtainable to represent all the series that generate that structure. This way, for each structure that we choose, we can start from a representative series and manipulate it with our operations for compositive criteria, yielding a new series that proliferates into a PP of the same structure. We will show the classification  for the particular, most interesting case of $n=12$ and $t=0$, showing the structures with their representative series and ordered by the amount of proliferations that they generate, that is to say, order of the PP, skipping those PPs of order 1 or 2. We will also show an example of the process of taking a representative series and manipulating it with the operations, observing that the structure of the PP does not change. Figure~\ref{fig:Classification} provides a base series from each equivalence class that can be used to obtain any other series in said class, as we have done in Figure~\ref{fig:Classification_use}.


In Figure~\ref{fig:Classification_use} we show a use case of the classification from Figure~\ref{fig:Classification}. We start from the last series of the classification, which generates a PP with the structure $[\{12\}]$, and then manipulate it to obtain the series used in "Accent in Pink", a piece written by Isabel Tardón and inspired by the painting of the same name by Wassily Kandinsky.

\section{Conclusion}

For each transformation and each transposition of any series of $n$ notes, we have found that there exists a pattern in the way that the proliferating process happens when we use the original series and its transformation as generators.

When the second generator series is a prime or an inversion, the process does not depend at all on the specific series chosen and it also does not generate any new material when compared to that which be can obtain through classical serialism. In the case of inversion, the order of the permutation is always 2, regardless of the transposition. 
For primes, the only transformation to go from the original series to the second generator is transposition. In this case, if $t$ is the number of tone-fractions transposed and $n$ is the number of notes in the series, then the number of proliferations obtained is always $n/gcd(n, t)$ and they are all primes with respect to the original series.

For retrogrades and retrograde inversions, the variety is far greater. These two cases do produce material which was not attainable by traditional serialism. 

In terms of composition, the variety of the series obtained in this way allow for large-scale musical projects to be constructed using one single dodecaphonic series as a base without the music becoming monotonous or uninteresting. For every series obtained by proliferation, there is a 12x12 matrix of all its transformations and, due to the method of construction, they will not all coincide, resulting in the material being multiplied. 

In the case of retrograde inversion, the parity in the number of notes and in the transposition both play a key role in determining the possible structures for the permutations: 

	\begin{enumerate}
	\item[(i)] When $n$ is even and $t$ is odd the PPs will all have a structure that is like a permutation of $n/2$ notes but repeating each cycle twice
	\item[(ii)]  When $n$ and $t$ are both even, all the PPs will have a cycle formed by an even number of notes $k$, and then two cycles for each cycle in an arbitrary permutation of $(n-k)/2$ notes. 
	\item[(iii)] When $n$ is odd, the behaviour is exactly the same as when $n$ and $t$ are both even, except that $k$ is now an odd number. If $t=0$, there is one additional restriction, which is that $k$ must be greater than 1.
\end{enumerate}

For retrograding only, we have several cases: 

    \begin{enumerate}
    \item[(i)] When the transposition applied is $t=0$, the permutation has order 2 and generates no new material. 
    
    \item[(ii)]When $n$ is an even number and $t=n/2$, the process behaves exactly like RI with the same number of notes $n$ and an odd transposition $t$. 

    \item[(iii)]When use R and a transposition $t$ such that $gcd(n, t)=1$, all the PPs obtainable are the ones where the number of cycles does not exceed $(n/2)+1$ and the number of cycles with an even number of notes is even if $n$ is 1 or 2 modulo 4, and it is odd when $n$ is 3 or 0 modulo 4.
    \item[(iv)]For the general case, where the transposition is any at all, the structures also follow, with little variation, the restrictions in (iii), but there are further restrictions that become too complex and varied to provide a general rule. We can provide, though, an inductive process to construct them.
\end{enumerate}
In pieces with shorter range, using a series and its proliferations creates music with a high level of internal coherence, but with the materials being related in a more complex way than usual, the listener has their expectations continuously shattered by the variety and the different kinds of sounds that can be achieved.

Overall, proliferating series are a natural evolution of the concept of serialism that we know. With this paper, we give composers the knowledge, tools and confidence to dive into this method of composition and explore its endless possibilities.

\section*{Acknowledgements}

The writing of this manuscript would not have been possible without the support of composer Reyes Oteo, who first introduced the concept of proliferating series to us, and the guidance of Cristina Draper and María de los Ángeles Gómez Molleda. We would also like to acknowledge Sebastián de la Torre Gil-Delgado for his suggestions on the code and our proofreaders Osiris García Parras, Irene Pastor González and José Luis Pérez Martín  who kindly gave their feedback and suggestions after reading the final drafts of the manuscript.

\bibliographystyle{apacite}
\bibliography{./interactapasample.bib}

\appendix
\section{Python code} 
\label{sec: appendix}
This code is also available on github at the following link with a more detailed explanation: \textcolor{blue}{\url{https://github.com/PabloMartinSantamaria/Proliferating-Series}}

The first code block automatizes the construction of proliferating series, and also has a function to collect the data of every possible proliferating series in a folder.

\begin{lstlisting} [style=PythonCustom]
from math import gcd
from functools import reduce
import os
from collections import defaultdict
from itertools import permutations
from contextlib import ExitStack
from enum import Enum
class Transformation(Enum):
    P, I, R, RI = "P", "I", "R", "RI"
    def apply(self, series, t, n):
        seq = reversed(series) if "R" in self.name else series
        return [((-x if "I" in self.name else x)+t)%n for x in seq]
def proliferatingPermutation(T, series, t, n):
    transformed = T.apply(series, t, n)
    idx = {value: i for (i, value) in enumerate(series)}
    return [transformed[idx[value]] for value in range(n)]
def lcm_list(series):
    return reduce(lambda a,b: a*b//gcd(a,b), series)
def cycleStructure(perm,n):
    visited = [False]*n
    cycles = []
    for i in range(n):
        if not visited[i]:
            current=perm[i]
            length=1
            while current!=i:
                visited[current]=True
                current=perm[current]
                length+=1
            cycles.append(length)
    return lcm_list(cycles), sorted(cycles)
def proliferations(T, series, t): 
    print("Proliferations of", series, "using", T.name,
          "and a transposition of", t, "tone-fractions:")
    print(series)
    n=len(series)
    perm = proliferatingPermutation(T, series, t, n)
    current = [perm[x] for x in series]
    while current!=series:
        print(current)
        current = [perm[x] for x in current]
    (order,cycles)=cycleStructure(perm,n) 
    print("Order:", order)
    print("Structure:", cycles)
def proliferations_data(T, t, n, path):
    path0 = os.path.join(path, T.name, f"Proliferations_{n}_notes")
    dirs = {"Data_Structures":os.path.join(path0,"Data_Structures"),
        "CompleteList":os.path.join(path0,"CompleteList"),
        "Data_Orders":os.path.join(path0,"Data_Orders")}
    for d in dirs.values():
        os.makedirs(d, exist_ok=True)
    paths = {key: os.path.join(dirs[key], f"transposition{t}.txt") 
             for key in dirs}
    occurrence_orders = defaultdict(int)
    occurrence_structures = defaultdict(int)
    with ExitStack() as stack:
        files = {key: stack.enter_context(open(paths[key], "w")) 
                 for key in paths}
        for seriesWithout0 in permutations(range(1, n)):
            series = (0,) + (seriesWithout0)
            perm = proliferatingPermutation(T, series, t, n)
            order, structure = cycleStructure(perm, n)
            files["CompleteList"].write(f"{series} --> {order}\n")
            occurrence_orders[order] += n
            occurrence_structures[(order, tuple(structure))] += n
        for k, v in sorted(occurrence_orders.items()):
            files["Data_Orders"].write(f"{k}: {v}\n")
        for k, v in sorted(occurrence_structures.items()):
            files["Data_Structures"].write(f"{k}: {v}\n")
\end{lstlisting}
Example 1: executing this cell will show the proliferations obtained when we apply a transformation of our choice to an arbitrary series. The list of numbers represents the series and the rightmost number represents the transposition.
\begin{lstlisting} [style=PythonCustom]
proliferations(Transformation.RI,[0,5,4,1,11,10,3,8,2,7,9,6],7)
proliferations(Transformation.R,[0,6,8,5,7,11,4,3,9,10,1,2],2)
\end{lstlisting}
Example 2: creates a folder in "desiredPath" with data of all possible proliferating permutations up to desired number of notes. Running the program for up to 12 notes will last several hours.
\begin{lstlisting} [style=PythonCustom]
for n in range(1,9):
    for t in range(n):
        proliferations_data(Transformation.RI,t,n,r"desiredPath")
        proliferations_data(Transformation.R ,t,n,r"desiredPath")
        proliferations_data(Transformation.P ,t,n,r"desiredPath")
        proliferations_data(Transformation.I ,t,n,r"desiredPath")
\end{lstlisting}
Code for calculating the list of all the possible proliferating permutations fixed a number of notes and a transformation:
\begin{lstlisting} [style=PythonCustom]
from math import gcd
from functools import reduce
from itertools import permutations
from collections import defaultdict
from copy import deepcopy
def lcm(series):
    return reduce(lambda a,b: a*b//gcd(a,b), series)
def d(l):
    return [i for i in l for _ in range(2)]
def canonicalList(L):
    return sorted(list(set([(lcm(l),tuple(l)) for l in L])))
def part(n, max=None):
    if max is None:       max = n
    if n == 0:            return [[]]
    if n < 0 or max == 0: return []
    withMax = part(n - max, max)
    withoutMax = part(n, max - 1)
    return [p + [max] for p in withMax] + withoutMax
def possiblePermutationsRI(n,t):
    L=[]
    if n==1: return [(1,[1])]
    elif n%2==0:
        if t%2==0: 
            for i in range(2,n+1,2): 
                L.extend([sorted([i]+d(l)) for l in part((n-i)//2)])
        else: L.extend([d(l) for l in part(n//2)])
    else: 
        for i in range(3,n+1,2): 
            L.extend([sorted([i]+d(l)) for l in part((n-i)//2)])
        if t%n!=0: L.extend([[1]+d(l) for l in part((n-1)//2)])
    return canonicalList(L)
def possiblePermutationsR_coprime_t(n):
    def evens(l):       return len([i for i in l if i%2==0])%2 
    def sgn(n):         return 1 if n%4 in [0,3] else 0
    L=[l for l in part(n) if len(l)<=(n/2)+1 and evens(l)==sgn(n)]
    return canonicalList(L)
def equalCycles(cycle1, cycle2): 
    if len(cycle1) != len(cycle2): return False
    sep = "|"
    s1 = sep.join(map(str, cycle1))
    s2 = sep.join(map(str, cycle2))
    return s1 in (s2 + sep + s2)
def GTPermutations(GT): 
    seps = [i for i in range(1,len(GT)) if GT[i]!=GT[i-1]]+[len(GT)]
    L=list(permutations(range(seps[0])))
    for i in range(1,len(seps)):
        newBlock=list(permutations(range(seps[i-1],seps[i])))
        L = [x+y for x in L for y in newBlock]
    return L 
def auxS(S1,S2,n,GT): 
    for perm in GTPermutations(GT):
        comprobation = [False for _ in range(len(S1))]
        flag=False 
        if n%2==1: 
            (len1,cyc1)=S1[0]
            (len2,cyc2)=S2[0]
            if len1!=len2: return False 
            else:
                cyc2 = [perm[i] for i in cyc2] 
                if cyc2==cyc1: comprobation[0]=True 
                else: flag=True
        for i in range(n%2,len(S1)):
            if flag: break
            (len1,cyc1)=S1[i]
            (len2,cyc2)=S2[i]
            if len1!=len2: return False
            else:
                cyc2 = [perm[i] for i in cyc2]
                if equalCycles(cyc2,cyc1): comprobation[i]=True
                else: flag=True
        if all(comprobation): return True
    return False   
def structPerm(S,n): 
    def cond(S,i): return S[i-1][0]!=S[i][0]
    if n%2==1: 
        seps=[1]+[i for i in range(2,len(S)) if cond(S,i)]+[len(S)]
    else:      
        seps=    [i for i in range(1,len(S)) if cond(S,i)]+[len(S)]
    L=list(permutations(S[:seps[0]]))
    for i in range(1,len(seps)):
        newBlock=list(permutations(S[seps[i-1]:seps[i]]))
        L = [x+y for x in L for y in newBlock]
    return [list(permutation) for permutation in L]
def equalStructures(S1,S2,n,GT):
    if len(S1)!=len(S2):                                         
        return False
    elif sorted([a[0] for a in S1])!=sorted([a[0] for a in S2]): 
        return False 
    else: 
        return any([auxS(S1,newS,n,GT) for newS in structPerm(S2,n)]) 
def cycleStructure(S):
    return tuple([length for (length,cycle) in S])
def eliminateRepeated(l,n,GT): 
    withoutRepeated = [] 
    dic=defaultdict(list) 
    for S in l:
        if n%2==1: sortedS=[S[0]]+sorted(S[1:]) 
        else:      sortedS=sorted(S)
        projection=cycleStructure(sortedS)
        if dic[projection]==[]:
            dic[projection]=[sortedS]
            withoutRepeated.append(sortedS)
        else: 
            flag=False
            for visited in dic[projection]:
                if equalStructures(visited,sortedS,n,GT): 
                    flag = True
                    break
            if not flag: 
                dic[projection].append(sortedS)
                withoutRepeated.append(sortedS)
    return withoutRepeated   
def auxR(n,GT): 
    if n==1: L = [[[1,[0]]]]
    elif n%2==1: 
        L=[]
        for i in range(len(GT)):
            if (i==len(GT)-1 or GT[i]!=GT[i+1]) and GT[i]==1: 
                auxGT = GT[:]
                del auxGT[i]
                for S in auxR(n-1,auxGT): 
                    L+=[[[1,[len(GT)-1]]]+S] 
            elif i==len(GT)-1 or GT[i]!=GT[i+1]:
                auxGT = GT[:]
                auxGT[i]-=1
                for S in auxR(n-1,auxGT): 
                    for j in range(len(S)):
                        (length,cycle)=S[j]
                        for k in range(length):
                            if cycle[k]==i: 
                                SNew=S[:]
                                SNew.pop(j)
                                newCycle=cycle[k:]+cycle[:k]+[i]
                                SNew.insert(0,[length+1,newCycle])
                                L+=[SNew]    
    else: 
        L=[]
        i = len(GT)-1 
        m = GT[i]
        if m == 1: 
            newGT=GT[:]
            del newGT[i]
            for S in auxR(n-1,newGT):
                SNew=S[:]
                SNew[0][0]+=1 
                SNew[0][1]+=[i]
                L+=[SNew]
        else:  
            newGT=GT[:]
            newGT[i]-=1
            for S in auxR(n-1,newGT):
                (l,c)=S.pop(0) 
                for j in range(0,len(S)): 
                    (l2,c2)=S[j]
                    for k in range(len(c2)):
                        if c2[k]==i:
                            SNew=S[:]
                            del SNew[j]
                            SNew.append([l+l2+1,[i]+c+c2[k:]+c2[:k]])
                            L+=[SNew]
                for j in range(len(c)): 
                    if c[j]==i:
                        SNew=S[:]
                        SNew.append([1+j,[i]+c[0:j]])
                        SNew.append([l-j,c[j:len(c)]])
                        L.append(SNew)            
    finalList=eliminateRepeated(L,n,GT) 
    return finalList
def possiblePermutationsR(n,GT):
    if isinstance(GT,int): GT=[n//gcd(n,GT) for _ in range(gcd(n,GT))]  
    L=[tuple(sorted([l for l,c in S])) for S in auxR(n,GT)]
    return canonicalList(L)
\end{lstlisting}
Examples of use: the first input is the number of notes and the second (if applicable) one is the transposition.  Functions for the cases of P and I are not included for.
\begin{lstlisting} [style=PythonCustom]
possiblePermutationsRI(12,5)
possiblePermutationsR_coprime_t(12)
possiblePermutationsR(12,3)
\end{lstlisting}

\end{document}